\documentclass{ws-p8-50x6-00}
\usepackage{amsmath}
\usepackage{amssymb}
\usepackage{amsfonts}

\newtheorem{thm}{Theorem}[section]

\newtheorem{defn}{Definition}[section]

\newfont{\w}{msbm9 scaled\magstep1}
\def\R{\mbox{\w R}}

\newcommand{\W}{\mathcal{W}}
\newcommand{\K}{\mathcal{K}}
\newcommand{\He}{\mathcal{H}}
\newcommand{\AK}{\mathcal{AK}}
\newcommand{\pd}{\partial}
\newcommand{\ddx}[1]{\frac{\pd}{\pd x^{#1}}}
\newcommand{\norm}[1]{||#1||}
\newcommand{\normq}[1]{\lVert#1\rVert ^2}
\newcommand{\Span}{{\rm span}}
\begin{document}

\title{SOME FOUR-DIMENSIONAL 
ALMOST HYPERCOMPLEX 
PSEUDO-HERMITIAN MANIFOLDS}

\author{MANCHO MANEV}
\address{Faculty of Mathematics and Informatics, University of Plovdiv \\
236 Bulgaria Blvd., Plovdiv 4003, Bulgaria\\
mmanev@pu.acad.bg,  mmanev@yahoo.com}
\author{KOUEI SEKIGAWA}
\address{Department of Mathematics, Faculty of Science\\
     Niigata University, Niigata, 950-2181, Japan\\
            sekigawa@sc.niigata-u.ac.jp}

\maketitle

\abstracts{ In this paper, a lot of examples of four-dimensional
manifolds with an almost hypercomplex pseudo-Hermitian structure
are constructed in several explicit ways. The received 4-manifolds
are characterized by their linear invariants in the known
aspects.}

\section*{Introduction}
In the study of almost hypercomplex manifolds the Hermitian
metrics are well known. The parallel study of almost hypercomplex
manifolds with skew-Hermitian metrics is in progress of
development~\cite{GrMaDi1},~\cite{GrMaDi2}.

Let $(M,H)$ be an almost hypercomplex manifold, i.e. $M$ is a
$4n$-dimen\-sional differentiable manifold and $H$ is a triple
$(J_1,J_2,J_3)$ of anticommuting almost complex structures, where
$J_3=J_1 \circ J_2$ \cite{So},\cite{AlMa}.

A standard hypercomplex structure for all $x(x^i,y^i,u^i,v^i)\in
T_pM$, $p\in M$ is defined in \cite{So} as follows
\begin{equation}\label{14}
\begin{array} {l}
J_1x(-y^i,x^i,v^i,-u^i),\quad J_2x(-u^i,-v^i,x^i,y^i),\quad
J_3x(v^i,-u^i,y^i,-x^i).
\end{array}
\end{equation}

Let us equip $(M,H)$ with a pseudo-Riemannian metric $g$ of
signature $(2n,2n)$ so that
\begin{equation}\label{1}
g(\cdot,\cdot)=g(J_1\cdot,J_1\cdot)=-g(J_2\cdot,J_2\cdot)=-g(J_3\cdot,J_3\cdot).
\end{equation}
We called such metric a \emph{pseudo-Hermitian metric} on an
almost hypercomplex manifold \cite{GrMaDi1}. It generates a
K\"ahler 2-form $\Phi$ and two pseudo-Hermitian metrics $g_2$ and
$g_3$ by the following way
\begin{equation}\label{G}
\Phi:=g(J_1\cdot,\cdot),\qquad g_2:=g(J_2\cdot,\cdot),\qquad
g_3:=g(J_3\cdot,\cdot).
\end{equation}
The metric $g$ ($g_2$, $g_3$, respectively) has an Hermitian
compatibility with respect to $J_1$ ($J_3$, $J_2$, respectively)
and a skew-Hermitian compatibility with respect to $J_2$ and $J_3$
($J_1$ and $J_2$, $J_1$ and $J_3$, respectively).

On the other hand, a quaternionic inner product $<\cdot,\cdot>$ in
$\mathbb{H}$ generates in a natural way the bilinear forms $g$,
$\Phi$, $g_2$ and $g_3$ by the following decomposition:
$<\cdot,\cdot>=-g+i\Phi+jg_2+kg_3$.

The structure $(H,G):=(J_1,J_2,J_3;g,\Phi,g_2,g_3)$ is called a
\emph{hypercomplex pseu\-do-Hermit\-ian structure} on $M^{4n}$ or
shortly a \emph{$(H,G)$-structure} on $M^{4n}$. The manifold
$(M,H,G)$ is called an \emph{almost hypercomplex pseudo-Hermitian
manifold} or shortly an \emph{almost
$(H,G)$-manifold}~\cite{GrMaDi1}.

The basic purpose of the recent paper is to construct explicit
examples of the $(H,G)$-manifolds of the lowest dimension at $n=1$
and to characterize them.

The following structural $(0,3)$-tensors  play basic role for the
characterization of the almost $(H,G)$-ma\-nifold
\[
F_\alpha (x,y,z)=g\bigl( \left( \nabla_x J_\alpha
\right)y,z\bigr)=\bigl(\nabla_x g_\alpha\bigr) \left( y,z \right),
\quad \alpha=1,2,3,
\]
where $\nabla$ is the Levi-Civita connection generated by $g$.

It is well known \cite{AlMa}, that the almost hypercomplex
structure $H=(J_\alpha)$ is a hypercomplex structure if the
Nijenhuis tensors
\[
N_\alpha(X,Y)=
    \left[X,Y \right]
    +J_\alpha\left[X,J_\alpha Y \right]
    +J_\alpha\left[J_\alpha X,Y \right]
    -\left[J_\alpha X,J_\alpha Y \right]
\]
vanish for each $\alpha=1,2,3$. Moreover, one $H$ is hypercomplex
iff two of $N_\alpha$ are zero.

Since $g$ is a Hermitian metric with respect to $J_1$, we use the
classification of the almost Hermitian manifolds given in
\cite{GrHe}. According to it the basic class of these manifolds of
dimension 4 are the class of almost K\"ahler manifolds $\AK=\W_2$
and the class of Hermitian manifolds $\He=\W_4$. The class of the
$\AK$-manifolds are defined by condition $d \Phi = 0$ or
equivalently $\mathop{\makebox{\huge $\sigma$}} \limits_{x,y,z}
F_1(x,y,z)=0$. The class of the Hermitian 4-manifolds is
determined by $N_1=0$ or
\[
\begin{array}{rl}
F_1(x,y,z)=\frac{1}{2} &\left[
g(x,y)\theta_1(z)-g(x,z)\theta_1(y)\right.\\[3pt]
&\left. -g(x,J_1y)\theta_1(J_1z)+g(x,J_1z)\theta_1(J_1y) \right]
\end{array}
\]
where $\theta_1(\cdot)=g^{ij}F_1(e_i,e_j,\cdot)=\delta\Phi(\cdot)$
for any basis $\{e_i\}_{i=1}^{4}$, and $\delta$ -- the
co\-de\-ri\-vative.

On other side, the metric $g$ is a skew-Hermitian one with respect
to $J_2$ and $J_3$. A classification of all almost complex
manifolds with skew-Hermitian metric (Norden metric or B-metric)
is given in \cite{GaBo}. The basic classes are:
\[
\begin{array}{l}
\begin{array}{rl}
\mathcal{W}_1: F_\alpha(x,y,z)=\frac{1}{4} &
    \left[
    g(x,y)\theta_\alpha(z)
    +g(x,z)\theta_\alpha(y)\right.\\[3pt]
    &\left.
    +g(x,J_\alpha y)\theta_\alpha(J_\alpha z)
    +g(x,J_\alpha z)\theta_\alpha(J_\alpha y)\right],\\[3pt]
\end{array}
\\[3pt]
\mathcal{W}_2: \mathop{\makebox{\huge $\sigma$}} \limits_{x,y,z}
F_\alpha(x,y,J_\alpha z)=0,\qquad
\mathcal{W}_3: \mathop{\makebox{\huge $\sigma$}} \limits_{x,y,z}
F_\alpha(x,y,z)=0,
\end{array}
\]
where \(\theta_\alpha(\cdot)=g^{ij}F_\alpha (e_i,e_j,\cdot)\),
\(\alpha =2,3\), for an arbitrary basis \(\{e_i\}_{i=1}^{4}\).

We denote the main subclasses of the respective complex manifolds
by $\W(J_\alpha)$, where $\W(J_1):=\W_4(J_1)$ \cite{GrHe}, and
$\W(J_\alpha):=\W_1(J_\alpha)$ for $\alpha=2,3$ \cite{GaBo}.

In the end of this section we recall some known facts from
\cite{GrMaDi1} and \cite{GrMaDi2}.

A sufficient condition an almost $(H,G)$-manifold to be an
integrable one is following

\begin{thm} \label{t31}
Let $(M,H,G)$ belongs to $ \W(J_\alpha) \bigcap \W(J_\beta)$. Then
$(M,H,G)$ is of class $\W(J_\gamma)$ for all cyclic permutations
$(\alpha, \beta, \gamma)$ of $(1,2,3)$.
\end{thm}

A pseudo-Hermitian manifold is called a
\emph{pseudo-hyper-K\"ahler manifold} (denotation $(M,H,G) \in
\K$), if $F_\alpha=0$ for every $\alpha=1,2,3$, i.e. the manifold
is K\"ahlerian with respect to each $J_\alpha$ (denotation
$(M,H,G) \in \K(J_\alpha)$).

\begin{thm}\label{t33}
If $(M,H,G) \in \K(J_\alpha)\bigcap \W(J_\beta)$ $(\alpha\neq\beta
\in \{1,2,3\})$ then $(M,H,G) \in \K$ .
\end{thm}
As $g$ is an indefinite metric, there exists isotropic vector
fields $X$ on $M$. Following \cite{GRMa} we consider the
invariants
\[
\norm{\nabla J_\alpha}^2= g^{ij}g^{kl}g\bigl( \left( \nabla_{e_i}
J_\alpha \right) e_k, \left( \nabla_{e_j} J_\alpha \right) e_l
\bigr), \qquad \alpha=1,2,3,
\]
where $\{e_i\}_{i=1}^{4}$ is an arbitrary basis of $T_pM$, $p\in
M$.

\begin{defn}
An $(H,G)$-manifold is called: (i) isotropic K\"ahlerian with
respect to $J_\alpha$
    if $\norm{\nabla J_\alpha}^2=0$ for some $\alpha\in\{1,2,3\}$;
(ii) isotropic hyper-K\"ahlerian
    if it is isotropic K\"ahlerian with respect to every $J_\alpha$ of $H$.
\end{defn}

\begin{thm}\label{t34}
Let $M$ be an $(H,G)$-manifold of class $\W=\bigcap_\alpha
\W(J_\alpha)$ $(\alpha=1,2,3)$ and $\norm{\nabla J_\alpha}^2$
vanishes for some $\alpha=1,2,3$. Then $(M,H,G)$ is an isotropic
hyper-K\"ahler manifold, but it is not pseudo-hyper-K\"ahlerian in
general.
\end{thm}

A geometric characteristic of the pseudo-hyper-K\"ahler manifolds
according to the curvature tensor $R=[\nabla ,\nabla ] -
\nabla_{[\ ,\ ]}$ induced by the Levi-Civita connection is given
in \cite{GrMaDi2}.
\begin{thm}\label{R=0}
Each pseudo-hyper-K\"ahler manifold is a flat pseudo-Riemann\-ian
manifold with signature $(2n,2n)$.
\end{thm}


\section{The two known examples of almost $(H,G)$-manifolds}

\subsection{A pseudo-Riemannian spherical manifold with
$(H,G)$-struc\-ture}

Following \cite{Wo} we have considered in \cite{GrMaDi1} and
\cite{GrMaDi2} a pseudo-Riemannian spherical manifold $S_2^4$ in
pseudo-Euclid\-ean vector space $\R_2^5$ of type $(--+++)$. The
structure $H$ is introduced on $\tilde S_2^4=S_2^4\setminus
\{(0,0,0,0,\pm 1)\}$ as in (\ref{14}) and the pseudo-Riemannian
metric $g$ is the restriction of the inner product of $\R_2^5 $ on
$\tilde{S}_2^4$. Therefore $\tilde{S}_2^4$ admits an almost
hypercomplex pseudo-Hermitian structure. The corresponding
manifold is of the class $\W(J_1)$ but it does not belong to $\W$
and it has a constant sectional curvature $k=1$.
Moreover, we established that the considered manifold is
conformally equivalent to a flat $\K(J_1)$-manifold, which is not
a
$\K$-manifold and $(\tilde{S}_2^4,H,G)$ is an Einstein manifold.


\subsection{The Thurston manifold with $(H,G)$-structure}

In \cite{GrMaDi1} we have followed the interpretation of Abbena
\cite{Ab} of the Thurston manifold. We have considered a
4-dimensional compact homogenous space $L\slash\Gamma$, where $L$
is a connected Lie group and $\Gamma$ is the discrete subgroup of
$L$ consisting of all matrices whose entries are integers. We have
introduced the almost hypercomplex structure
$H=\left(J_\alpha\right)$ on $T_EL$ as in (\ref{14}) and we
translate it on $T_AL$, $A\in L$, by the action of the left
invariant vector fields. The $J_\alpha$ are invariant under the
action of $\Gamma$, too. By analogy we have defined a left
invariant pseudo-Riemannian inner product in $T_EL$. It generates
a pseudo-Riemannian metric $g$ on $M^4=L$. Then the generated
4-manifold $M$ is equipped with a suitable $(H,G)$-structure and
$(M,H,G)$ is a $\W (J_1)$-manifold but it does not belong to the
class $\W$.


\section{Engel manifolds with almost $(H,G)$-structure}

In the next two examples we consider
$M=\mathbb{R}^4=\left\{(x^1,x^2,x^3,x^4)\right\}$ with a basis \(
\begin{array}{c}
\left\{e_1=\ddx{1},e_2=\ddx{2}+x^1\ddx{3}+x^3\ddx{4},e_3=-\ddx{3},e_4=-\ddx{4}\right\}
\end{array}
\) and an Engel structure $\mathcal{D}=\Span \{e_1,e_2\}$, i.e. an
absolutely non-integrable regular two-dimensional distribution on
$TM$ \cite{GRMa}.

\subsection{Double isotropic hyper-K\"ahlerian structures but
neither hypercomplex nor symplectic}

At first we use the introduced there a pseudo-Riemannian metric
and almost complex structures given by
\begin{equation}\label{41}
\begin{array}{c}
g=(dx^1)^2+\{1-(x^1)^2-(x^3)^2\}(dx^2)^2-(dx^3)^2   \\[3pt]
-(dx^4)^2-2x^1dx^2dx^3+2x^3dx^2dx^4,
\\[3pt]
J:\; Je_1=e_2,\;Je_2=-e_1,\;Je_3=e_4,\;Je_4=-e_3,
\\[3pt]
J':\; J'e_1=e_2,\;J'e_2=-e_1,\;J'e_3=-e_4,\;J'e_4=e_3.
\end{array}
\end{equation}

It is given in \cite{GRMa} that $(J,g)$ and $(J',g)$ are a pair of
indefinite almost Hermitian structures which are isotropic
K\"ahler but neither complex nor symplectic.

It is clear that \( \left\{e_i\right\}_{i=1}^4 \) is an
orthonormal $(++--)$-basis. We accomplish the introduction of an
$(H,G)$-structure on $M$ by
\[
J_1:=J';\;
J_2:\;J_2e_1=e_3,\;J_2e_2=e_4,\;J_2e_3=-e_1,\;J_2e_4=-e_2;\;
J_3:=J_1J_2.
\]

By direct computations we verify that the constructed manifold is
an $(H,G)$-manifold and it is isotropic hyper-K\"ahlerian but not
K\"ahlerian and not integrable  with non-vanishing Lie forms with
respect to any $J_\alpha\ (\alpha=1,2,3)$.

\begin{rem}
If we define $J_1$ as $J$ instead of $J'$ then the kind of example
is not changed. So we receive a pair of almost $(H,G)$-structures
corresponding to the given almost complex structures.
\end{rem}

The non-zero components of the curvature tensor $R$ and the basic
linear invariant of the almost Hermitian manifold $(M,J_1,g)$ are
given in \cite{GRMa} by
\[
\begin{array}{c}
R_{1221}=\frac{3}{4},\;
R_{1331}=-R_{2142}=-R_{2442}=-R_{3143}=R_{3443}=\frac{1}{4},\;
R_{2332}=1;
\\[3pt]
\normq{F_1}=0,\quad \normq{N_1}=8,\quad\tau=0,\quad\tau^*_1=-2,
\end{array}
\]
where the following denotations are used for
$\varepsilon_a=\normq{e_a}$
\[
\begin{array}{c}
\normq{F_1}=\normq{\nabla\Phi}=\sum_{a,b,c=1}^4 \varepsilon_a
\varepsilon_b \varepsilon_c F_1(e_a,e_b,e_c)^2,\\[3pt]
\normq{N_1}=\sum_{a,b=1}^4 \varepsilon_a \varepsilon_b
\normq{N_1(e_a,e_b)},\\[3pt]
 \tau=\sum_{a,b=1}^4 \varepsilon_a
\varepsilon_b R(e_a,e_b,e_b,e_a),\quad
\tau^*_1=\frac{1}{2}\sum_{a,b=1}^4 \varepsilon_a \varepsilon_b
R(e_a,J_1e_a,e_b,J_1e_b).
\end{array}
\]
We get the corresponding linear invariants with respect to $J_2$
and $J_3$:
\[
\begin{array}{lll}
&\normq{F_2}=0,\qquad \normq{N_2}=0,\qquad&\tau^*_2=0;
\\
&\normq{F_3}=0,\qquad \normq{N_3}=-8,\quad&\tau^*_3=0,
\end{array}
\]
where $\tau^*_\alpha=\sum_{a,b=1}^4 \varepsilon_a \varepsilon_b
R(e_a,e_b,J_\alpha e_b,e_a);\;\alpha=2,3$.


\subsection{Double isotropic hyper-K\"ahlerian structures which
are non-integrable but symplectic}

Now we consider the same Engel manifold
$(M=\mathbb{R}^4,\mathcal{D})$ but let the pseudo-Riemannian
metric and the pair of almost complex structures be defined by
other way: \cite{GRMa}
\[
\begin{array}{c}
g=(dx^1)^2-\{1-(x^1)^2+(x^3)^2\}(dx^2)^2+(dx^3)^2\\[3pt]
-(dx^4)^2-2x^1dx^2dx^3+2x^3dx^2dx^4,
\\[3pt]
J:\; Je_1=e_3,\;Je_2=e_4,\;Je_3=-e_1,\;Je_4=-e_2,
\\[3pt]
J':\; J'e_1=e_3,\;J'e_2=-e_4,\;J'e_3=-e_1,\;J'e_4=e_2.
\end{array}
\]

In this case \( \left\{e_i\right\}_{i=1}^4 \) is an orthonormal
basis of type $(+-+-)$. It is shown that $(M,J,g)$ and $(M,J',g)$
are a pair of isotropic K\"ahler almost K\"ahler manifolds with
vanishing linear invariants.

We accomplish the introduced almost complex structures to almost
hypercomplex structures on $M$ by using the following way: we set
the given $J$ (resp. $J'$) as $J_1$ (resp. $J'_1$), then we
introduce $J_2$ (resp. $J'_2$) by
\begin{equation}\label{42h}
\begin{array}{c}
J_2:\;J_2e_1=e_2,\;J_2e_2=-e_1,\;J_2e_3=-e_4,\;J_2e_4=e_3;\\[3pt]
J'_2:\;J'_2e_1=e_2,\;J'_2e_2=-e_1,\;J'_2e_3=e_4,\;J'_2e_4=-e_3
\end{array}
\end{equation}
and finally we set $J_3:=J_1J_2$ (resp. $J'_3:=J'_1J'_2$).

 It is easy to check that $H=(J_\alpha)$ and $H'=(J'_\alpha)$
 together with $g$ generate a pair of almost hypercomplex
 pseudo-Hermitian structures on $M$.

 We characterize the both received $(H,G)$-manifolds as isotropic
 hyper-K\"ahler but not K\"ahler manifolds
 and not integrable manifolds with non-vanishing Lie forms
 with respect to any $J_\alpha$.
 Moreover, they have the following linear invariants:
\[
\begin{array}{c}
\normq{N_1}=0,\;\; \normq{N_2}=-\normq{N_3}=8,\;\;
\normq{F_\alpha}=0,\;\; \tau=\tau^*_\alpha=0 \;\; (\alpha=1,2,3).
\end{array}
\]

\section{Real spaces with almost $(H,G)$-structure}
\subsection{Real semi-space with almost $(H,G)$-structure}

Let us consider the real semi-space
$\mathbb{R}_+^4=\left\{\left(x^1,x^2,x^3,x^4\right),\; x^i\in
\mathbb{R},\; x^1>0 \right\}$ with the basis given by
\(
\begin{array}{c}
\left\{e_1=x^1\ddx{1},e_2=x^1\ddx{2},e_3=x^1\ddx{3},e_4=x^1\ddx{4}\right\}.
\end{array}
\)
It is clear that this basis is orthonormal of type $(++--)$ with
respect to the pseudo-Riemannian metric \(
g=\left\{(dx^1)^2+(dx^2)^2-(dx^3)^2-(dx^4)^2\right\}/(x^1)^2. \)
We introduce an almost hypercomplex structure $H=(J_\alpha)$ as
follows
\begin{equation}\label{43h}
\begin{array}{l}
J_1:\;J_1e_1=e_2,\;J_1e_2=-e_1,\;J_1e_3=e_4,\;J_1e_4=-e_3;\\[3pt]
J_2:\;J_2e_1=e_3,\;J_2e_2=-e_4,\;J_2e_3=-e_1,\;J_2e_4=e_2;\quad J_3=J_1J_2\\[3pt]
\end{array}
\end{equation}
and we check that $H$ and $g$ generates an almost
$(H,G)$-structure on $\mathbb{R}_+^4$.

We verify immediately that $H$ is integrable and the obtained
hypercomplex pseudo-Hermit\-ian manifold $(\mathbb{R}_+^4,H,G)$
belongs to the class $\W=\bigcap_\alpha \W(J_\alpha)$  but it is
not isotropic K\"ahlerian with respect to $J_\alpha$
$(\alpha=1,2,3)$.

By direct computations we obtain for the curvature tensor that
$R=-\pi_1$, i.e. the manifold has constant sectional curvatures
$k=-1$ and it is an Einstein manifold. Moreover, the linear
invariants are
\[
\begin{array}{c}
\normq{N_\alpha}=0,\quad
2\normq{F_1}=4\normq{\theta_1}=-\normq{F_\beta}=-\normq{\theta_\beta}=16,\\[3pt]
\tau=-3\tau^*_1=-12,\; \tau^*_\beta=0,\;
\end{array}
\]
where $\alpha=1,2,3;\ \beta=2,3$; and $(\mathbb{R}_+^4,H,G)$ is
conformally equivalent to a pseudo-hyper-K\"ahler manifold by the
change $\bar{g}=(x^1)^2g$.

\subsection{Real quarter-space with almost $(H,G)$-structure}

Let the real quarter-space
\[
M=\mathbb{R}_+^2\times
\mathbb{R}_-^2=\left\{\left(x^1,x^2,x^3,x^4\right),\; x^i\in
\mathbb{R},\; x^1>0,\; x^3>0 \right\}\] be equipped with a
pseudo-Riemannian metric
\[
\begin{array}{c}
g=\frac{1}{(x^1)^2}\left\{(dx^1)^2+(dx^2)^2\right\}
-\frac{1}{(x^3)^2}\left\{(dx^3)^2+(dx^4)^2\right\}.
\end{array}
\]
Then the basis \(
\left\{e_1=x^1\ddx{1},e_2=x^1\ddx{2},e_3=x^3\ddx{3},e_4=x^3\ddx{4}\right\}
\) is an orthonormal one of type $(++--)$. We introduce an almost
hypercomplex structure $H=(J_\alpha)$ $(\alpha=1,2,3)$ as in the
previous example by (\ref{43h}).

The received almost $(H,G)$-manifold is a
$\mathcal{K}(J_1)$-manifold and an isotropic hyper-K\"ahler
manifold. As a corollary, $N_1=0$, $F_1=0$, $\theta_1=0$ and hence
$\normq{N_1}=\normq{F_1}=\normq{\theta_1}=0$. For the $J_\alpha$
$(\alpha =2,3)$ the Nijenhuis tensors $N_\alpha$, the tensors
$F_\alpha$, and the Lie forms $\theta_\alpha$ are non-zero
(therefore $H$ is not integrable), but the linear invariants
$\normq{N_\alpha}$, $\normq{F_\alpha}$  and
$\normq{\theta_\alpha}$ vanish.

The non-zero components of the curvature tensor are given by
$R_{1221}=-R_{3443}=-1$. For the Ricci tensor we have
$\rho_{ii}=-1$ $(i=1,...,4)$. Therefore the basic non-zero
sectional curvatures are $k(e_1,e_2)=-k(e_3,e_4)=-1$ and the
scalar curvatures $\tau$, $\tau^*_\alpha$ $(\alpha =1,2,3)$ are
zero.

\section{Real pseudo-hyper-cylinder with almost
$(H,G)$-structure}

Let $\mathbb{R}^5_2$ be a pseudo-Euclidean real space with an
inner
product $\langle\cdot,\cdot\rangle$ of signature 
$(+++--)$. Let us consider a pseudo-hyper-cylinder defined by
\[
S:\; (z^2)^2+(z^3)^2-(z^4)^2-(z^5)^2=1,
\]
where $Z\left(z^1,z^2,z^3,z^4,z^5\right)$ is the positional vector
at $p\in S$. We use the following parametrization of $S$ in the
local coordinates $\left(u^1,u^2,u^3,u^4\right)$ of $p$:
\[
\begin{array}{c}
Z=Z(u^1,\ \cosh u^4 \cos u^2,\ \cosh u^4 \sin u^2,\ \sinh u^4 \cos
u^3,\ \sinh u^4 \sin u^3).
\end{array}
\]

We consider a manifold on the surface
$\tilde{S}=S\setminus\{u^4=0\}$. Then the basis
\( \left\{e_1=\pd_{1},\ e_2=\frac{1}{\cosh u^4}\pd_{2},\
e_3=\frac{1}{\sinh u^4}\pd_{3},\ e_4=\pd_{4}\right\} \)
of $T_p\tilde{S}$ at $p\in \tilde{S}$ is an orthonormal basis of
type $(++--)$ with respect to the restriction $g$ of
$\langle\cdot,\cdot\rangle$ on $\tilde{S}$. Here and further
$\pd_{i}$ denotes $\frac{\pd Z}{\pd u^i}$ for $i=1,...,4$;

We introduce an almost hypercomplex structure by the following way
\begin{equation}\label{45h}
\begin{array}{l}
J_1:\;J_1e_1=e_2,\;J_1e_2=-e_1,\;J_1e_3=-e_4,\;J_1e_4=e_3;\\[3pt]
J_2:\;J_2e_1=e_3,\;J_2e_2=e_4,\;J_2e_3=-e_1,\;J_2e_4=-e_2;\quad J_3=J_1J_2\\[3pt]
\end{array}
\end{equation}
and check that $H=(J_\alpha)$ and the pseudo-Riemannian metric $g$
generate an almost $(H,G)$-structure on $\tilde{S}$.

By straightforward calculations with respect to $\{e_i\}$
$(i=1,...,4)$ we receive that the almost $(H,G)$-manifold
$\tilde{S}$ is not integrable with non-zero Lie forms regarding
any $J_\alpha$ of $H$ and it has the following linear invariants:
\[
\begin{array}{l}
\normq{N_1}=2\normq{F_1}=2\normq{\nabla
J_1}=8\normq{\theta_1}=-8\tanh^2 u^4;\\[3pt]
\normq{N_2}=-8\coth^2 u^4,\quad \normq{\theta_2}=\left(2\tanh
u^4+\coth u^4\right)^2,\\[3pt]
\normq{F_2}=\normq{\nabla J_2}=4\left(2\tanh^2 u^4+\coth^2
u^4\right);\\[3pt]
\normq{N_3}=-8\left(\tanh u^4-\coth u^4\right)^2,\quad
\normq{\theta_3}=\left(\tanh
u^4+\coth u^4\right)^2,\\[3pt]
\normq{F_3}=\normq{\nabla J_3}=4\left(\tanh^2 u^4+\coth^2
u^4\right).
\end{array}
\]

The non-zero components of the curvature tensor and the
corresponding Ricci tensor and scalar curvatures are given by
\[
\begin{array}{c}
R_{2332}=-1,\quad R_{2442}=-\tanh^2 u^4,\quad R_{3443}=\coth^2
u^4\\[3pt]
\rho_{22}=1+\tanh^2 u^4,\quad \rho_{33}=-1-\coth^2 u^4,\quad
\rho_{44}=-\tanh^2 u^4-\coth^2 u^4\\[3pt]
\tau=2\left(1+\tanh^2 u^4+\coth^2 u^4\right),\quad
\tau^*_\alpha=0,\quad \alpha=1,2,3.
\end{array}
\]

Hence $(\tilde{S},H,G)$ has zero associated scalar curvatures and
$H$ is a non-integrable structure on it.
%

\section{Complex surfaces with almost $(H,G)$-structure}

The following three examples concern several surfaces
$S^2_\mathbb{C}$ in a 3-dimensional complex Euclidean space
$\left(\mathbb{C}^3,\langle\cdot,\cdot\rangle\right)$. It is well
known that the decomplexification of $\mathbb{C}^3$ to
$\mathbb{R}^6$ using the $i$-splitting, i.e.
$(Z^1,Z^2,Z^3)\in\mathbb{C}^3$, where $Z^k=x^k+iy^k$ $(x^k,y^k\in
\mathbb{R})$, is identified with $(x^1,x^2,x^3,y^1,y^2,y^3)\in
\mathbb{R}^6$. Then the multiplying by $i$ in $\mathbb{C}^3$
induces the standard complex structure $J_0$ in $\mathbb{R}^6$.
The real and the opposite imaginary parts of the complex Euclidean
inner product $\Re\langle\cdot,\cdot\rangle$ and
$-\Im\langle\cdot,\cdot\rangle$ are the standard skew-Hermitian
metrics $g_0$ and $\tilde{g}_0=g_0(\cdot,J_0\cdot)$ in
$(\mathbb{R}^6,J_0,g_0,\tilde{g}_0)$, respectively. So, the
natural decomplexification of an $n$-dimensional complex Euclidean
space is the $2n$-dimensional real space with a complex
skew-Hermitian structure $(J_0,g_0,\tilde{g}_0)$.

\subsection{Complex cylinder with almost $(H,G)$-structure}

Let $S^2_\mathbb{C}$ be the cylinder in
$\left(\mathbb{C}^3,\langle\cdot,\cdot\rangle\right)$ defined by
$(Z^1)^2+(Z^2)^2=1$. Let the corresponding surface $S^4$ in
$(\mathbb{R}^6,J_0,g_0,\tilde{g}_0)$ be parameterized as follows
\[
\begin{array}{rl}
S^4:\; Z=Z(&\cos u^1 \cosh u^3,\ \sin u^1 \cosh u^3,\ u^2,\\[3pt]
&\sin u^1 \sinh u^3,\ -\cos u^1 \sinh u^3,\ u^4).
\end{array}
\]
Then the local basis $\left\{\pd_1,...,\pd_4\right\}$ is
orthonormal of type $(++--)$ and it generates the metric
$g=(du^1)^2+(du^2)^2-(du^3)^2-(du^4)^2$ on $S^4$. The almost
hypercomplex structure $H$ is determined as usually by (\ref{14}).
It is easy to verify that the received $(H,G)$-manifold is a flat
pseudo-hyper-K\"ahler manifold.

\subsection{Complex cone with almost $(H,G)$-structure}

Now let $S^2_\mathbb{C}$ be the complex cone in
$\left(\mathbb{C}^3,\langle\cdot,\cdot\rangle\right)$ determined
by the equation
$(Z^1)^2+(Z^2)^2-(Z^3)^2=0$ 
. Then we consider the corresponding 4-dimensional surface $S$ in
$(\mathbb{R}^6,J_0,g_0,\tilde{g}_0)$ by the following
parametrization of $Z$:
\[
\begin{array}{l}
(u^1 \cos u^2 \cosh u^4-u^3 \sin u^2 \sinh u^4,\ u^1 \sin
u^2 \cosh u^4+u^3 \cos u^2 \sinh u^4,\ u^1,\\[3pt]
\;u^1 \sin u^2 \sinh u^4+u^3 \cos u^2 \cosh u^4,-u^1 \cos u^2
\sinh u^4+u^3 \sin u^2 \cosh u^4,\ u^3).
\end{array}
\]
Further we consider a manifold on
$\tilde{S}=S\setminus\{0,0,0,0,0,0\}$, i.e. we exclude the plane
$u^1=u^3=0$ from the domain of $S$ which maps the origin. Then the
derived metric $g$ on $\tilde{S}$ has the following non-zero
components regarding $\left\{\pd_k\right\}$:
\[
g_{11}=-g_{33}=2,\quad g_{22}=-g_{44}=(u^1)^2-(u^3)^2,\quad
g_{24}=g_{42}=2u^1u^3.
\]
We receive the following orthonormal basis of signature $(++--)$:
\[
\begin{array}{c}
\left\{e_1=\frac{1}{\sqrt{2}}\pd_{1},\
e_2=\lambda\pd_{2}+\mu\pd_{4},\ e_3=\frac{1}{\sqrt{2}}\pd_{3},\
e_4=-\mu\pd_{2}+\lambda\pd_{4}\right\},
\end{array}
\]
where 
$\lambda=u^1/\{(u^1)^2+(u^3)^2\}$, $\mu=u^3/\{(u^1)^2+(u^3)^2\}$.
We introduce a structure $H$ as in (\ref{14}). It is easy to check
that $H$ and $g$ generate an almost $(H,G)$-structure on
$\tilde{S}$. By direct computations we get that the received
$(H,G)$-manifold is a flat hypercomplex manifold which is
K\"ahlerian with respect to $J_1$ but it does not belong to
$\mathcal{W}(J_2)$ or $\mathcal{W}(J_3)$ and the Lie forms
$\theta_2$ and $\theta_3$ are non-zero. The corresponding linear
invariants are given by
\[
\begin{array}{c}
\normq{F_2}=\normq{\nabla
J_2}=2\normq{\theta_2}=16\left\{\mu^2-\lambda^2\right\},\\[3pt]
\normq{F_3}=\normq{\nabla
J_3}=2\normq{\theta_3}=4\left\{\mu^2-\lambda^2\right\}.
\end{array}
\]

\subsection{Complex sphere with almost $(H,G)$-structure}

In this case let $S^2_\mathbb{C}$ be the unit sphere in
$\left(\mathbb{C}^3,\langle\cdot,\cdot\rangle\right)$ defined by
$(Z^1)^2+(Z^2)^2+(Z^3)^2=1$. After that we consider the
corresponding 4-surface $S$ in
$(\mathbb{R}^6,J_0,g_0,\tilde{g}_0)$ with the following
parametrization of 
$Z(x^1,x^2,x^3,y^1,y^2,y^3)$:
\[S:\;
\begin{array}{l}
x^1=\cos u^1 \cos u^2 \cosh u^3 \cosh u^4-\sin u^1 \sin u^2 \sinh u^3 \sinh u^4, \\[3pt]
x^2=\cos u^1 \sin u^2 \cosh u^3 \cosh u^4+\sin u^1 \cos u^2 \sinh u^3 \sinh u^4, \\[3pt]
x^3=\sin u^1 \cosh u^3, \\[3pt]
y^1=\cos u^1 \sin u^2 \cosh u^3 \sinh u^4+\sin u^1 \cos u^2 \sinh u^3 \cosh u^4, \\[3pt]
y^2=-\cos u^1 \cos u^2 \cosh u^3 \sinh u^4+\sin u^1 \sin u^2 \sinh u^3 \cosh u^4, \\[3pt]
y^3=-\cos u^1 \sinh u^3.
\end{array}
\]
Further we consider a manifold on $\tilde{S}=S\setminus\{0,0,\pm
1,0,0,0\}$, i.e. we exclude the set $u^1=\pm\pi/2,u^3=0$ from the
domain $(-\pi,\pi)^2\times\mathbb{R}^2$ of $S$ which maps the pair
of "poles".

The induced metric on $\tilde{S}$ has the following non-zero local
components:
\[
\begin{array}{c}
g_{11}=-g_{33}=1,\quad g_{22}=-g_{44}=\cos^2 u^1 \cosh^2
u^3-\sin^2 u^1 \sinh^2
u^3,\\[3pt]
g_{24}=g_{42}=2\sin u^1\cos u^1 \sinh u^3\cosh u^3.
\end{array}
\]
Further we use the following orthonormal basis of signature
$(++--)$:
\[
\begin{array}{c}
\left\{e_1=\pd_{1},\ e_2=\lambda\pd_{2}+\mu\pd_{4},\ e_3=\pd_{3},\
e_4=-\mu\pd_{2}+\lambda\pd_{4}\right\},
\end{array}
\]
where 
%
$\lambda=\frac{\cos u^1 \cosh u^3}{\cos^2 u^1+\sinh^2 u^3}$,
$\mu=\frac{\sin u^1 \sinh u^3}{\cos^2 u^1+\sinh^2 u^3}$.
We introduce a structure $H$ as in (\ref{14}) and we verify that
$H$ and $g$ generate an almost $(H,G)$-structure on $\tilde{S}$.
By direct computations we get that $(\tilde{S},H,G)$ is a
$\mathcal{K}(J_2)$-manifold of pointwise constant totally real
sectional curvatures
\[
\begin{array}{c}
\nu=\frac{\sinh^2 2u^3-\sin^2 2u^1}{4(\cos^2 u^1+\sinh^2
u^3)^4},\quad \nu^*_2=\frac{\sin 2u^1\sinh 2u^3}{2(\cos^2
u^1+\sinh^2 u^3)^4},
\end{array}
\] where
$\nu:=\frac{R(x,y,y,x)}{\pi_1(x,y,y,x)}$,
$\nu^*_2:=\frac{R(x,y,y,J_2x)}{\pi_1(x,y,y,x)}$ for a basis
$\{x,y\}$ of any non-de\-ge\-ne\-rate totally real section
$\sigma$ (i.e. $\sigma\perp J_2\sigma$). $(\tilde{S},J_2,g)$ is an
almost Einstein manifold since its Ricci tensor is $\rho=2(\nu g-
\nu^*_2 g_2)$. But, the Nijenhuis tensors and the Lie forms
corresponding to other two almost complex structures $J_1$ and
$J_3$ are non-zero. Beside that, we receive the following linear
invariants:
\[
\begin{array}{c}
\tau=8\nu,\qquad \tau^*_1=0,\qquad \tau^*_2=8\nu^*_2,\qquad \tau^*_3=0,\\[3pt]
\normq{N_1}=2\normq{\nabla J_1}=8\normq{\theta_1}=-32\nu,\;
-\normq{N_3}=2\normq{\nabla J_3}=8\normq{\theta_3}=32\nu.
\end{array}
\]

\section{Lie groups with almost $(H,G)$-structure}

The next two examples are inspired from an example of a locally
flat almost Hermitian surface constructed in \cite{TrVa}. Let
$\mathcal{L}$ be a connected Lie subgroup of
$\mathcal{GL}(4,\mathbb{R})$ consisting of matrices with the
following non-zero entries
\[
\begin{array}{c}
a_{11}=a_{22}=\cos u_{1},\quad a_{12}=-a_{21}=\sin
u_{1},\\[3pt]
a_{13}=u_{2},\quad a_{23}=u_{3},\quad a_{33}=1,\quad a_{44}=\exp
u_{4}
\end{array}
\]
for arbitrary $u^{1},u^{2},u^{3},u^{4}\in \mathbb{R}$.

The Lie algebra of $\mathcal{L}$ is isomorphic to the Lie
subalgebra of $\mathfrak{gl}(4;\mathbb{R})$ generated by the
matrices $X_1$, $X_2$, $X_3$, $X_4$ with the the following
non-zero entries:
\[
(X_1)_{13}=(X_2)_{12}=-(X_2)_{21}=(X_3)_{23}=(X_4)_{44}=1.
\]

\subsection{A Lie group as a complex manifold but non-hypercomplex one}
For the first recent example let us substitute the following
pseudo-Riemannian $g$ for the metric on $\mathcal{L}$ used there:
$g(X_i,X_j)=\varepsilon_a \delta_{ij}$, where $1\leq i,j\leq 4$;
$\varepsilon_1=\varepsilon_2=-\varepsilon_3=-\varepsilon_4=1$.
Further we introduce an $\mathcal{L}$-invariant almost
hypercomplex structure $H$ on $\mathcal{L}$ as in (\ref{14}).
Then, there is generated an almost $(H,G)$-structure on
$\mathcal{L}$ and the received manifold is complex with respect to
$J_2$ but non-hypercomplex and the Lie forms do not vanish. The
non-zero components of the curvature tensor $R$ is determined by
$R_{1221}=R_{1331}=-R_{2332}=1$ and the linear invariants are the
following:
\[
\begin{array}{c}
\normq{N_1}=2\normq{\nabla J_1}=8\normq{\theta_1}=-\normq{\nabla
J_2}=-2\normq{\theta_2}=\normq{N_3}=-8,\\[3pt]
\normq{\nabla J_3}=12\normq{\theta_3}=12,\qquad
\tau=-\tau^*_1=2,\quad \tau^*_2=\tau^*_3=0.
\end{array}
\]
\subsection{A Lie group as a flat K\"ahler manifold but non-hypercomplex one}
For the second example we use the following pseudo-Riemannian $g$
on $\mathcal{L}$: $g(X_i,X_j)=\varepsilon_a \delta_{ij}$, where
$1\leq i,j\leq 4$;
$\varepsilon_1=-\varepsilon_2=\varepsilon_3=-\varepsilon_4=1$. We
actually substitute only the type of the signature: $(+-+-)$ for
$(++--)$ of the basis $\{X_1,X_2,X_3,X_4\}$. Then we introduce $H$
by the following different way:
\[
\begin{array}{llll}
J_1X_1=X_3,\; &J_1X_2=X_4,\; &J_1X_3=-X_1,\; &J_1X_4=-X_2,\\[3pt]
J_2X_1=-X_4,\; &J_2X_2=X_3,\; &J_2X_3=-X_2,\; &J_2X_4=X_1,\qquad
J_3=J_1J_2.
\end{array}
\]
Therefore we obtain that the constructed $(H,G)$-manifold is flat
and K\"ahler\-ian with respect to $J_1$ but regarding $J_2$ and
$J_3$ it is not complex and the structural tensors have the form
$F_2(X,Y,Z)=-\theta_2(J_3X)g(Y,J_3Z)$,
$F_3(X,Y,Z)=-\theta_3(J_2)g(Y,J_2Z)$. The non-zero linear
invariants for $\beta=2,3$ are the following: \(
\begin{array}{c}
-\normq{N_\beta}=2\normq{\nabla
J_\beta}=2\normq{F_\beta}=8\normq{\theta_\beta}=8.
\end{array}
\)

\section*{Acknowledgement} The first author is supported by
the Matsumae International Foundation's fellowship under the
guidance of the second author.

\end{document}